\documentclass{amsart} 
\usepackage{graphicx} 
\newcommand{\eps}{\varepsilon} 
 
\def\C{\mathbb C} 
\def\N{\mathbb N} 
\def\dimH{\operatorname{dim_H}} 
\def\dimP{\operatorname{dim_P}}

\def\dim{\operatorname{dim}} 
\def\diam{\operatorname{diam}} 
\newtheorem{thm}{Theorem} 
\theoremstyle{remark} 
\newtheorem*{rem}{Remark} 
\newtheorem*{ack}{Acknowledgment} 
\begin{document} 
\title[The set where the iterates are bounded]{On the set where the iterates of 
an entire function are bounded} 
\author{Walter Bergweiler}\thanks{Supported 
by a Chinese Academy of Sciences Visiting 
Professorship for Senior International Scientists, Grant 
No.\ 2010 TIJ10. Also supported by  
the Deutsche Forschungsgemeinschaft, Be~1508/7-1, 
the EU Research Training Network CODY and 
the ESF Networking Programme HCAA.}  
\email{bergweiler@math.uni-kiel.de} 
\address{Mathematisches Seminar, 
Christian--Albrechts--Universit\"at zu Kiel, 
Lude\-wig--Meyn--Str.~4, 
D--24098 Kiel, 
Germany} 
\begin{abstract} 
We show that for a transcendental entire function the set 
of points whose orbit under iteration is 
bounded can have arbitrarily small positive
Hausdorff dimension. 
\end{abstract} 
\maketitle 
\section{Introduction} 
The main objects studied in complex dynamics are the \emph{Fatou set} 
$F(f)$ of a rational or entire function $f$,  defined as the set of 
all points where the iterates $f^n$ of $f$ form a normal family, and 
the \emph{Julia set} $J(f)$, which  is the complement of $F(f)$. 
 In the 
dynamics of transcendental entire functions -- and this is the case 
we shall be concerned with -- a fundamental role is also played by 
the escaping set 
$$ 
I(f)=\{z\in\C: f^n(z)\to\infty  \text{ as }n\to\infty\}. 
$$ 
 The first systematic study of this set was undertaken 
 by Eremenko~\cite{Eremenko89} who, among other results, showed 
 that $I(f)\neq\emptyset$ and in fact $I(f)\cap J(f)\neq\emptyset$ 
 for every transcendental entire function~$f$.  Moreover, $J(f)=\partial 
 I(f)$. In this paper we will consider the set 
$$ 
K(f)=\{z\in\C: (f^n(z)) \text{ is bounded}\}. 
$$ 
As repelling periodic points are dense in the Julia set~\cite{Baker68}, the 
properties of $I(f)$ mentioned above also hold for $K(f)$; that is, 
$K(f)\cap J(f)\neq\emptyset$ and $J(f)=\partial K(f)$. 
 
For a polynomial $f$ the set $K(f)$ is called the \emph{filled Julia set} 
of~$f$ and we have $K(f)=\C\setminus I(f)$, but for a transcendental entire 
function $f$ there are points which are neither in $K(f)$ nor in 
$I(f)$, for example there are points in $J(f)$ whose orbit is dense in $J(f)$. 
However, there may also be points in $F(f)$ which are neither in 
$K(f)$ nor in $I(f)$; see~\cite[Example~1]{Eremenko87}. 
 
We denote the Hausdorff dimension and the packing dimension of a 
subset $A$ of $\C$ by $\dimH A$ and $\dimP A$, respectively. We 
refer to 
Falconer's book~\cite{Falconer90} 
for the definition 
of these dimensions and further information. Here we only note that  
we always  have 
$\dimH A\leq\dimP A$; 
 see~\cite[p.~48]{Falconer90}. By a result of Baker~\cite{Baker75}, the 
Julia set of a transcendental entire function $f$ contains continua. 
 In fact, even $I(f)\cap J(f)$ contains continua  and thus $\dim_H (I(f)\cap J(f))\geq 1$; 
 cf.~\cite[Theorem~5]{Rippon05a} and~\cite[Theorem~1.3]{Rippon10}.

 A major 
open question in transcendental dynamics is whether 
 $\dim_H J(f)>1$ 
for every transcendental entire function~$f$. It was proved by 
Stallard (\cite{Stallard96}, see also~\cite{BRS,Rippon06}) that this is the 
 case for functions in the Eremenko-Lyubich class $B$ which consists of 
all transcendental entire functions for which the set of critical 
values and finite asymptotic values is bounded. Bara\~nski, 
Karpi\~nska and Zdunik~\cite{Baranski09} 
showed that for $f\in B$ there exists a compact, invariant 
Cantor subset $C$ of $J(f)$ with $\dimH C>1$. 
In particular, $\dimH (K(f)\cap J(f))>1$ 
for $f\in B$. 
 On the other hand, Rempe and Stallard~\cite{Rempe10} showed that 
 there are functions $f\in B$ for which $\dimH I(f)=1$. 
 
We consider the dimensions of $K(f)$ for entire functions 
which need not be in Eremenko-Lyubich class. 
The following result is a special case of a result  
of Rempe~\cite[Corollary~2.11]{Rempe09} who proved 
that the hyperbolic dimension of an Ahlfors islands map 
is positive. 
\begin{thm}\label{thm1} 
If $f$ is a transcendental entire function, then 
$\dimH (K(f)\cap J(f))>0$. 
\end{thm} 
Theorem~\ref{thm1} is also implicit in Stallard's~\cite{Stallard94}  
proof that $\dimH J(f)>0$ for transcendental meromorphic functions~$f$. 
The proofs in~\cite{Rempe09,Stallard94} 
 are based on suitable versions of the Ahlfors islands theorem; 
see~\cite[Theorem~6.2]{Hayman64} or, for an alternative 
proof,~\cite{Bergweiler98}.  This is used to to obtain an iterated 
function scheme (see~\cite{Falconer90}), 
 whose invariant set is a (hyperbolic) Cantor subset  
of $K(f)\cap J(f)$ which can be shown to have positive Hausdorff 
dimension. We note, however, that for entire and meromorphic  
functions different versions of the  Ahlfors islands theorem 
have to be used; see the discussion in~\cite[Section~6.4]{Bergweiler00}. 
For entire functions such a hyperbolic, invariant Cantor subset  
of $K(f)\cap J(f)$ is also constructed in~\cite{Christensen96}. 
 
It is the purpose of this note to show that Theorem~\ref{thm1} is 
best possible even for entire functions. 
\begin{thm}\label{thm2} 
For every $\eps>0$ there exists a transcendental entire function $f$ 
such that 
$\dimH K(f)\leq \dim_P K(f)<\eps.$ 
\end{thm} 
For an introduction to the dynamics of transcendental entire (and 
meromorphic) functions we refer to~\cite{Bergweiler93}. Results on 
dimensions of Julia sets of transcendental functions 
are surveyed in~\cite{Stallard08}. 
 
\begin{ack} 
I thank Lasse Rempe, Phil Rippon and Gwyneth Stallard for helpful comments. 
\end{ack} 
 
\section{Proof of Theorem~\ref{thm2}}\label{proof2} 
Let $C$ be a large positive constant and define $(a_k)_{k\geq 1}$ 
recursively by $a_1=1$ and 
\begin{equation}\label{p2a} 
a_{k+1}=8 C \, a_k \prod_{j=1}^{k-1}\frac{a_k}{a_j} 
\end{equation} 
for $k\geq 1$. (Here $\prod_{j=1}^0 a_1/a_j =1$ so that $a_2=8 C 
a_1=8 C$.) Induction shows that $(a_k)$ increases and that 
\begin{equation}\label{p2b} 
\frac{a_{k+1}}{a_k}\geq 8 C 
\prod_{j=1}^{k-1}\frac{a_k}{a_{k-1}}\geq (8 C)^k 
\end{equation} 
for all $k$. Thus 
$$ 
f(z)= C \, z \prod_{k=1}^\infty\left(1-\frac{z}{a_k}\right) 
$$ 
defines an entire function~$f$. For $k\geq 1$ we put 
$$ 
r_k=\frac{2k+1}{2k+2}a_k 
 \quad\text{and}\quad 
 s_k=10 a_k 
$$ 
and we set $r_0=0$ and $s_0=16/C$. 
For large $C$ we have $r_k<s_k<r_{k+1}$ for $k\geq 0$.
We define, for $k\geq 0$,
$$ 
A_k=\{z\in\C: r_k\leq |z|\leq s_{k}\} 
 \quad\text{and}\quad 
 B_k=\{z\in\C: s_k< |z|< r_{k+1}\}. 
$$ 
We will show that 
\begin{equation}\label{p2f} 
f(B_k)\subset B_{k+1} 
\end{equation} 
for all $k\geq 1$. In order to do so we note first that 
by~\eqref{p2b} we can achieve that 
\begin{equation}\label{p2g} 
\frac{a_{k+1}}{a_k}> 320 e (k+1)\geq 2k+4 
\end{equation} 
for all $k\geq 1$ by choosing $C$ sufficiently large. 
We deduce that if $1\leq j\leq k-1$, then 
$(2k+2)a_j\leq (2k+2)a_{k-1}\leq  a_k$ and hence 
\begin{equation}\label{p2i} 
1+\frac{r_k}{a_j}\leq 
\frac{a_k}{(2k+2)a_j}+\frac{r_k}{a_j}=\frac{a_k}{a_j} 
\end{equation} 
and 
\begin{equation}\label{p2h} 
\frac{r_k}{a_j}-1\geq \frac{r_k}{a_j}-\frac{a_k}{(2k+2)a_j}= 
\frac{k}{k+1}\frac{a_k}{a_j}. 
\end{equation} 
Moreover, it follows from~\eqref{p2b} that we can achieve that
\begin{equation}\label{p2j} 
 \prod_{j=k+1}^\infty\left(1+\frac{10a_k}{a_j}\right)\leq 2 
 \quad\text{and}\quad 
  \prod_{j=k+1}^\infty\left(1-\frac{10a_k}{a_j}\right)\geq\frac{9}{10}\geq 
  \frac12 
\end{equation} 
for all $k\geq 1$ by choosing $C$ large. 
 
 For $k\geq 1$ we deduce from~\eqref{p2a}, \eqref{p2i} 
 and~\eqref{p2j} that if $|z|=r_k$, then 
$$ 
\begin{aligned} 
|f(z)| & \leq 
 C\,  r_k \prod_{j=1}^{k-1}\left(1+\frac{r_k}{a_j}\right) 
\cdot \left(1+\frac{r_k}{a_k}\right) \cdot 
\prod_{j=k+1}^\infty\left(1+\frac{r_k}{a_j}\right) 
  \\ 
& \leq 4 C\, a_k \prod_{j=1}^{k-1}\frac{a_k}{a_j} 
 =\frac12 a_{k+1} 
 < r_{k+1}. 
\end{aligned} 
$$ 
Similarly, \eqref{p2a}, \eqref{p2g}, \eqref{p2h} 
 and~\eqref{p2j} yield that if $|z|=r_k$, then 
\begin{equation}\label{p2l} 
\begin{aligned} 
|f(z)| & \geq 
  C \, r_k \prod_{j=1}^{k-1}\left(\frac{r_k}{a_j}-1\right) 
\cdot \left(1-\frac{r_k}{a_k}\right) \cdot 
\prod_{j=k+1}^\infty\left(1-\frac{r_k}{a_j}\right)  \\ 
&\geq C 
\left(\frac{k}{k+1}\right)^{k} a_k \prod_{j=1}^{k-1}\frac{a_k}{a_j} 
  \cdot \frac{1}{2k+2}\cdot \frac12\\ 
&\geq \frac{C}{2e(2k+2)}  a_k \prod_{j=1}^{k-1}\frac{a_k}{a_j} 
= \frac{a_{k+1}}{32e (k+1)} 
> 10 a_{k} 
=s_k. 
\end{aligned} 
\end{equation} 
The last two inequalities imply that 
\begin{equation}\label{p2m} 
f(z)\in B_{k} \quad\text{for }|z|=r_k 
\end{equation} 
if $k\geq 1$.  Next we note that if $k\geq 1$ and and $|z|=s_k$, then 
\begin{equation}\label{p2n} 
\begin{aligned} 
|f(z)| & \geq 
 C \, s_k \prod_{j=1}^{k-1}\left(\frac{s_k}{a_j}-1\right) 
\cdot \left(\frac{s_k}{a_k}-1\right) \cdot 
\prod_{j=k+1}^\infty\left(1-\frac{s_k}{a_j}\right)  \\ 
&\geq  10 C\, a_k \prod_{j=1}^{k-1}\frac{9a_k}{a_j} \cdot 9 \cdot 
\frac{9}{10} 
= \frac{9^{k+1}}{8} a_{k+1} 
>s_{k+1}. 
\end{aligned} 
\end{equation} 
Similarly as in~\eqref{p2j} we also see that if $|z|=s_0=16/C$, then 
\begin{equation}\label{p2o} 
|f(z)|  \geq 
 C \, s_0 
\prod_{j=1}^\infty\left(1-\frac{s_0}{a_j}\right) 
 \geq  C \, s_0 
 \frac{9}{10} =\frac{16\cdot 9}{10}>10=s_1,
\end{equation} 
provided $C$ is chosen large enough.
Also, since $s_k< r_{k+1}$ for all $k\geq 0$, we deduce 
from~\eqref{p2m}, with $k$ replaced by $k+1$, that $|f(z)|< r_{k+2}$ 
for $|z|=s_k$. Together with~\eqref{p2n} and~\eqref{p2o} this yields 
that
\begin{equation}\label{p2p} 
f(z)\in B_{k+1} \quad\text{for }|z|=s_k 
\end{equation} 
if $k\geq 0$. 
Combining this  with~\eqref{p2m} we obtain~\eqref{p2f}. 
 
Next we show that with $L=C/(4e)$ we have 
\begin{equation}\label{p2q} 
|f'(z)|\geq  2^{k}L 
  \quad\text{for}\ z\in A_k. 
\end{equation} 
In order to do so we note first that if $p$ is a real polynomial 
with real zeros, then each interval bounded by two adjacent zeros 
 of $p$ contains exactly one zero of $p'$, and besides multiple 
 zeros of $p$ there are no further zeros of~$p'$. In particular, 
 $p'$ has only real zeros. Moreover, we see that $p$ has no positive 
 local minima and no negative local maxima. 
 
Since our function $f$ is a limit of real polynomials with real, 
non-negative zeros, $f'$ is also a limit of such polynomials. It 
follows  that $f'$ has no positive local minima and no negative 
local maxima. This implies that if a compact interval contains no 
zero of $f'$, then $|f'|$ assumes its minimum in the interval at one 
of the endpoints of the interval. 
 The fact that $f'$ is a limit of real polynomials 
with real, non-negative zeros also implies that  $|f'|$ takes its 
minimum on a circle around the origin at the intersection of this 
circle with the positive real axis. We will see that $f'$ has no 
zeros in the intervals $[r_k,s_k]$. The above arguments then imply 
that 
\begin{equation}\label{p2r} 
\min_{z\in A_k}|f'(z)|=\min \{|f'(r_k)|,|f'(s_k)|\}. 
\end{equation} 
In order to prove that $f'$ has no zeros in the intervals 
$[r_k,s_k]$, we note that if $r_k\leq x<a_k$ and $1\leq j\leq k-1$, 
 then $x>2a_j$ by~\eqref{p2b} and hence $x/(x-a_j)<2$.
Thus
\begin{equation}\label{p2s} 
\begin{aligned} 
\frac{xf'(x)}{f(x)} 
 & =1 + \sum_{j=1}^\infty \frac{x}{x-a_j} 
 \leq 1 + \sum_{j=1}^{k-1} \frac{x}{x-a_j} + \frac{r_k}{r_k-a_k} \\ 
& \leq 1 + 2(k-1) -  (2k+1) =-2<0 \quad\text{for}\ r_k\leq x<a_k.
\end{aligned} 
\end{equation} 
On the other hand, using~\eqref{p2b} it is not difficult to see that 
by choosing $C$ large we can achieve that if $k\geq 1$, then 
\begin{equation}\label{p2t} 
\frac{xf'(x)}{f(x)} \geq 1 - \sum_{j=k+1}^\infty \frac{s_k}{a_j- 
s_k}\geq\frac{1}{2} \quad\text{for}\ a_k<x\leq s_k. 
\end{equation} 
With $a_0=0$ this also holds for $k=0$ if $C$ is large. It follows 
from~\eqref{p2s} and~\eqref{p2t} that $f'$ has no zeros in the 
intervals $[r_k,s_k]$. Thus~\eqref{p2r} holds. Moreover, 
\eqref{p2b}, \eqref{p2l} and~\eqref{p2s} yield that
\begin{equation}\label{p2u} 
|f'(r_k)|\geq 2 \frac{|f(r_k)|}{r_k}\geq 2 \frac{a_{k+1}}{32e(k+1)a_k} \geq 
\frac{(8C)^k}{16e(k+1)}\geq\frac{C}{4e}2^k= 2^k L 
\end{equation} 
for $k\geq 1$ while \eqref{p2b}, \eqref{p2n} and~\eqref{p2t} give 
\begin{equation}\label{p2v} 
|f'(s_k)|\geq \frac{1}{2} \frac{|f(s_k)|}{s_k}\geq \frac{1}{2} 
\frac{9^{k+1}a_{k+1}}{80a_k} 
\geq \frac{1}{2} \frac{9^{k+1}(8C)^k}{80} 
\geq 4C\,2^k \geq 2^k L 
\end{equation} 
for $k\geq 1$. Finally, $f'(0)=C\geq L$ and~\eqref{p2o} implies
that
\begin{equation}\label{p2w} 
|f'(s_0)|\geq \frac{1}{2} \frac{s_{1}}{s_0}= \frac{10}{32}C\geq L. 
\end{equation} 
Now~\eqref{p2q} follows from~\eqref{p2r}, \eqref{p2u}, 
 \eqref{p2v} and~\eqref{p2w}. 
 
To estimate the dimension of $K(f)$, we fix $N\in\N$ and put 
$$ 
K_N(f)=\{z\in\C:|f^n(z)|\leq s_{N}\text{ for all }n\in\N\} 
$$ 
It follows from~\eqref{p2f} that $K_N(f)$ consists of all points $z$ 
for which 
 $f^n(z)\in\bigcup_{k=0}^{N} A_k $ 
 for all $n\in\N$. 
Thus, assuming that $C$ is chosen such that $L=C/(4e)>1$,
we deduce from~\eqref{p2q} that $K_N(f)$ is a conformal repeller; 
see~\cite[Chapter~8]{Przytycki10} and~\cite[Chapter~5]{Zinsmeister00} for the 
definition and 
 properties of conformal repellers. 
 It follows (see \cite[Corollary~8.1.7]{Przytycki10} 
 or \cite[Theorem~5.12]{Zinsmeister00}) that the Minkowski dimension, 
packing dimension and Hausdorff dimension of $K_N(f)$  all coincide 
and are given by Bowen's formula. This formula says that with 
$F=f|_{K_N(f)}$ these dimensions are equal to the unique zero of the 
pressure function $t\to P(F,t)$ defined by 
$$ 
P(F,t)  =\lim_{n\to\infty} \frac{1}{n} \log 
 \left( \sum_{z\in F^{-n}(a)}|(F^n)'(z)|^{-t}\right), 
$$ 
for some $a\in K_N(f)$. 
 
In order to apply Bowen's formula we note that 
every point in $K_N(f)$ has $N+1$ preimages under~$F$. 
Let $a\in A_k$. 
It follows from~\eqref{p2m} and the maximum principle that $F$ 
has no $a$-points in $A_j$ for $0\leq j\leq k-2$. 
Moreover, it follows from~\eqref{p2m} and~\eqref{p2p} 
that $F$ and $F-a$ have the same number of zeros in $A_j$ 
for $k\leq j\leq N$. Thus $F$ has exactly one $a$-point in 
$A_j$ for $k\leq j\leq N$. 
We conclude that $a$ has $k-1$ preimages under $F$ 
in~$A_{k-1}$. 
It follows from the above discussion, 
 together with~\eqref{p2q}, that for $a\in K_N(f)$ and 
$t>0$ we have 
$$ 
\sum_{b\in F^{-1}(a)}|F'(b)|^{-t}
\leq \sum_{k=0}^N \left(2^kL\right)^{-t} \leq 
L^{-t}\sum_{k=0}^\infty 2^{-tk} =  \frac{L^{-t}}{1-2^{-t}}. 
$$ 
Now 
$$ 
\begin{aligned} 
 \sum_{z\in F^{-(n+1)}(a)}|(F^{n+1})'(z)|^{-t} 
& = \sum_{b\in F^{-1}(a)} \sum_{z\in 
F^{-n}(b)}|(F^{n+1})'(z)|^{-t}\\ 
& = \sum_{b\in F^{-1}(a)}|F'(b)|^{-t} 
 \sum_{z\in F^{-n}(b)} |(F^{n})'(z)|^{-t}.
 \end{aligned} 
$$ 
With 
$$ 
S_n(t)=\sup_{c\in K_N(f)} \sum_{z\in F^{-n}(c)}|(F^{n})'(z)|^{-t} 
$$ 
we thus have 
$$ 
S_{n+1}(t)\leq  \frac{L^{-t}}{1-2^{-t}} S_n(t). 
$$ 
Induction shows that 
\begin{equation}\label{p2a2} 
\sum_{z\in F^{-n}(a)}|(F^n)'(z)|^{-t}\leq S_n(t)\leq 
\left(\frac{L^{-t}}{1-2^{-t}}\right)^n 
\end{equation} 
for all $a\in K_N(f)$. Thus 
\begin{equation}\label{p2a3} 
P(F,t)\leq \log  \frac{L^{-t}}{1-2^{-t}}. 
\end{equation} 
Given $t>0$, we can achieve that the right hand side of~\eqref{p2a3} 
is negative by choosing $C$ and hence $L$ large. Then 
the Minkowski, packing and Hausdorff 
dimension of $K_N(f)$ are less than~$t$ for all~$N$. 
 Since 
$K(f)=\bigcup_{N=1}^\infty K_N(f)$, we deduce that $\dim_P K(f)\leq 
t$. As $t>0$ can be chosen arbitrarily small, the conclusion 
follows. 
 
\begin{rem} 
The thermodynamic formalism of~\cite{Przytycki10,Zinsmeister00} is not 
actually needed to obtain an \emph{upper}  bound for $\dimH K_N(f)$. 
 As $K_N(f)$ does not intersect the 
postcritical set of $F$, there exists $\delta>0$ such that Koebe's 
distortion theorem may be applied to all inverse branches of the 
iterates of $F$ on  
the disk $D(a,\delta)=\{z\in\C: |z-a|<\delta\}$. 
We obtain 
$$ 
F^{-n}(D(a,\delta))\subset \bigcup _{z\in F^{-n}(a)} 
 D\left(z,\frac{C}{|(F^n)'(z)|}\right) 
$$ 
for some constant~$C$. Now~\eqref{p2a2} shows that 
$F^{-n}(D(a,\delta))$ can be covered by $(N+1)^n$ sets $V_j$ whose 
diameters satisfy 
$$\sum_j (\diam V_j)^t\leq (2C)^t 
\left(\frac{L^{-t}}{1-2^{-t}}\right)^n.$$ The compact set $K_N(f)$ can 
be covered by finitely many, say~$M$, disks $D(a,\delta)$. Hence we 
obtain a covering of $K_N(f)=F^{-n}(K_N(f))$ by $M(N+1)^n$ sets 
$W_j$ satisfying 
$$\sum_j (\diam W_j)^t\leq M(2C)^t 
\left(\frac{L^{-t}}{1-2^{-t}}\right)^n.$$ 
\end{rem} 
This implies that the $t$-dimensional Hausdorff measure of $K_N(f)$ 
is~$0$, provided $L$ is again chosen such that $L^{-t}<1-2^{-t}$.

\end{document}